\documentclass{article}

\usepackage[english]{babel}
\usepackage{amssymb}
\usepackage{amsthm}

\newtheorem{theorem}{Theorem}[section]
\newtheorem{definition}[theorem]{Definition}
\newtheorem{lem}[theorem]{Lemma}
\newtheorem{prop}[theorem]{Proposition}
\newtheorem{corollary}[theorem]{Corollary}

\newtheorem{ex}[theorem]{Example}

\newcommand{\hh}{{\mathbb{H}}}

\newcommand{\rr}{{\mathbb{R}}}
\newcommand{\nn}{{\mathbb{N}}}

\newcommand{\s}{{\mathbb{S}}}

\newcommand{\z}{{\mathcal{Z}}}

\title{\bf The open mapping theorem for regular quaternionic functions}
\author{ Graziano Gentili\footnote{Partially supported by G.N.S.A.G.A.
of the I.N.D.A.M. and by M.I.U.R.}  \\ \normalsize Dipartimento di
Matematica ``U. Dini'', Universit\`a di Firenze \\ \normalsize
Viale Morgagni 67/A, 50134 - Firenze, Italy,  gentili@math.unifi.it \\
\and Caterina Stoppato \\ \normalsize Dipartimento di Matematica
``U. Dini'', Universit\`a di Firenze \\ \normalsize Viale Morgagni
67/A, 50134 - Firenze, Italy,  stoppato@math.unifi.it\\}

\date{  }

\begin{document}
\maketitle


\begin{abstract}
The basic results of a new theory of regular functions of a
quaternionic variable have been recently stated, following an
idea of Cullen. In this paper we prove the minimum modulus principle
and the open mapping theorem for regular functions. The proofs involve
some peculiar geometric properties of such functions which are of
independent interest.
\end{abstract}


\section{Introduction}

Let $\hh$ denote the skew field of real quaternions. Its elements are
of the form $q=x_0+ix_1+jx_2+kx_3$ where the $x_l$ are real, and $i$,
$j$, $k$, are imaginary units (i.e. their square equals $-1$) such
that $ij=-ji=k$, $jk=-kj=i$, and $ki=-ik=j.$ 

After Hamilton's discovery of quaternions, the richness of the
theory of holomorphic functions of one complex variable, along with
motivations from physics, aroused 
interest in a theory of quaternion valued functions
of a quaternionic variable. In fact several interesting theories have been 
introduced in the last century.
The most famous is  
the one due to Fueter, \cite{fueter 1, fueter 2}, the basic results of
which are accurately summarized in \cite{sudbery}. Recent work on
Fueter-regularity includes \cite{libro shapiro} and reference therein.

A different theory of
quaternion-valued functions of one quaternionic variable has been
recently proposed by G. Gentili and D. C. Struppa \cite{cras,
advances}. The theory is based on a definition of regularity for
quaternionic functions inspired by C. G. Cullen \cite{cullen}.
Several basic results of the theory are proven in \cite{advances}, including
the Cullen-regularity of quaternion power series and some nice
properties of their zeros. The study of the zero-sets has been
further developed in \cite{zeros}.

We shall now quickly review the definition of Cullen-regular function
and the basic properties of such a function. Denote by $\s$ the
two-dimensional sphere of quaternion imaginary units: $\s = \{q \in
\hh : q^2 =-1\}$. For all imaginary unit $I \in \s$, let $L_I = \rr +
I \rr$ be the complex line through $0, 1$ and $I$.

\begin{definition}\label{definition}
Let $\Omega$ be a domain in $\hh$ and let $f : \Omega \to \hh$ be a
real differentiable function. $f$ is said to be
\textnormal{Cullen-regular} if, for all $I \in \s$, the function
$\bar \partial_I f : \Omega \cap L_I \to \hh$ defined by
\begin{equation}
\bar \partial_I f (x+Iy) = \frac{1}{2} ( \frac{\partial}{\partial
x}+I\frac{\partial}{\partial y} ) f_I (x+Iy)
\end{equation}
vanishes identically.
\end{definition}

With the notations $\Omega_I = \Omega \cap L_I$ and $f_I =
f_{|_{\Omega_I}}$, we may refer to the vanishing of $\bar \partial_I
f$ saying that the restriction $f_I$ is holomorphic on $\Omega_I$.
From now on we will omit Cullen's name and refer to these functions
just as regular functions. As observed in \cite{advances}, a
quaternion power series $\sum_{n \in \nn} q^n a_n$ with $a_n \in \hh$
defines a regular function in its domain of convergence, which proves
to be an open ball $B(0,R) = \{q \in \hh : |q| <R\}$. In the same paper, it
is proven that

\begin{theorem} 
If $f : B = B(0,R) \to \hh$ is regular then there exist quaternions
$a_n \in \hh$ such that
\begin{equation}
f(q)=\sum_{n \in \nn} q^n a_n
\end{equation}
for all $q \in B$. In particular, $f \in C^{\infty}(B)$.
\end{theorem}

Many classical results in complex analysis are
extended in \cite{advances} to regular functions $f : B(0,R) \to \hh$: the identity
principle, the maximum modulus principle, the Cauchy representation
formula, the Liouville theorem, the Morera theorem and the Schwarz
lemma. 

In this paper we prove the open mapping theorem for
regular functions. We first define an operation on such functions which 
allows us to deduce a weak version of the minimum modulus principle 
(see theorem \ref{minimum}) from the maximum modulus principle proven in
\cite{advances}. We then study the behaviour of a regular function $f$
restricted to any 2-sphere $S=x+y\s$ ($x,y \in \rr, y\neq 0$) and find out
that $f_{|_{S}}$ is either constant or an affine map of $S$ onto a
2-sphere $b+\s c$ with $b, c\in \hh$. This peculiar geometric property, 
which is of independent interest, leads to the following version of the 
minimum modulus principle:

\begin{theorem}[Minimum Modulus Principle]
Let $f : B=B(0,R) \to \hh$ be a regular function. If $|f|$ has a
local minimum point $p\in B$ then either $f(p)=0$ or $f$ is
constant.
\end{theorem}

We define the degenerate set $D_{f}$ of $f$ as the union of all the 
2-spheres $S=x+y\s$  ($x,y \in \rr, y\neq 0$) on which 
$f$ is constant. If $f$ is not constant then $D_f$ turns out to have empty 
interior and by means of the above result we prove

\begin{theorem}[Open Mapping Theorem]
Let $f : B(0,R) \to \hh$ be a non-constant regular function and let
$D_{f}$ be its degenerate set. Then $f: B(0, R)\setminus \overline{D_{f}}\to
\hh$ is  
open.
\end{theorem}

\begin{corollary} Let $f : B(0,R) \to \hh$ be a regular function. 
If $f$ is not constant on any 2-sphere $x+y\s$ with $y\neq 0$,  then $f$ is  
open.
\end{corollary}

We give an example of regular function $f : B=B(0,R) \to \hh$ which is open
on $B \setminus \overline{D_{f}}$ but not on $B$. Finally, we show that an 
open subset $U\subseteq B(0,R)$ such that $U=\bigcup_{x+yI\in U} x+y\s$ 
always maps to an open set $f(U)$, even when $U$ intersects the 
degenerate set $D_f$. In particular $f(B(0,R))$ is always an open set.


\section{Preliminary results}

We will recall in this section some results and definitions which
appear 
in the study of the structure of the
zero-sets of regular functions presented in \cite{advances, zeros}.
We begin with a simmetry property of the zeros proven in
\cite{advances}. This result and the following are proven for polynomials in
\cite{shapiro}.

\begin{theorem}
Let $f : B = B(0,R) \to \hh$ be a regular function and let $x, y \in
\rr, y>0$ be such that $x^2 +y^2 <R^2$. If there exist distinct
imaginary units $I,J \in \s$ such that $f(x+yI) = f(x+yJ) = 0$, then
$f(x+yK) = 0$ for all $K \in \s$.
\end{theorem}

In other words, if $f$ has more than one zero on the 2-sphere $x+y\s$
then it vanishes identically on $x+y\s$. In \cite{zeros} the zero-set
is further characterized as follows.

\begin{theorem}
Let $f$ be a regular function on an open ball $B = B(0,R)$. If $f$ is
not identically zero then its zero-set consists of isolated points or
isolated 2-spheres of the form $x + y \s$, for $x,y \in \rr, y\neq 0$.
\end{theorem}

The study of the zeros presented in \cite{zeros} involves the
following operations on regular functions $f : B(0,R) \to \hh$.

\begin{definition}\label{multiplication}
Let $f, g$ be regular functions on an open ball $B = B(0,R)$ and let
$f(q) = \sum_{n \in \nn} q^n a_n, g(q) = \sum_{n \in \nn} q^n b_n$ be
their power series expansions. We define the \textnormal{regular
product} of $f$ and $g$ as the regular function $f*g : B \to \hh$
defined by
\begin{equation}
f*g(q) = \sum_{n \in \nn} q^n c_n, \ c_n = \sum_{k=0}^n a_k b_{n-k}.
\end{equation}
We define the \textnormal{regular conjugate} of $f$ as the function
$f^c : B \to \hh$, defined by the following formula: $f^c(q) =
\sum_{n \in \nn} q^n \bar a_n$. Finally, we define the
\textnormal{symmetrization} of $f$, as $f^s = f * f^c = f^c * f$. 
\end{definition}

Note that the regular multiplication is an associative,
non-commutative binary operation. Since no confusion can arise, we
may also write $f(q)*g(q)$ for $f*g(q)$. Also note that the
symmetrization $f^s(q) = \sum_{n \in \nn} q^n \sum_{k = 0}^n a_k \bar
a_{n-k}$ has real coefficients. The zero-sets of $f^c$ and $f^s$ are
characterized in \cite{zeros} as follows.

\begin{theorem}\label{conjugatezeros}
Let $f$ be a regular function on $B = B(0,R)$. For all $x,y \in \rr$
with $x+y\s \subseteq B$, the zeros of the regular conjugate $f^c$ on
$x+y\s$ are in one-to-one correspondence with those of $f$. Moreover,
the symmetrization $f^s$ vanishes exactly on the sets $x+y\s$ on
which $f$ has a zero.
\end{theorem}

Note that $x+y\s$ is a 2-sphere if $y \neq 0$ and a real singleton
$\{x\}$ if $y= 0$.


\section{The reciprocal function}

We define a new operation on regular functions, associating to any 
regular $f : B(0,R) \to \hh$ a function $f^{-*}$ called its regular 
reciprocal. The reason for this name will become apparent later. 
Denote by $\z_f$ the zero set of a regular function $f$.

\begin{definition}
Let $f : B = B(0,R) \to \hh$ be a regular function and let $f^c,f^s$
be its regular conjugate and symmetrization. We call the
\textnormal{regular reciprocal} of $f$ the function $f^{-*} : B
\setminus \z_{f^s} \to \hh$ defined by
\begin{equation}
f^{-*}(q) = \frac{1}{f^s(q)} f^c(q).
\end{equation}
\end{definition}

We will use the shorthand notation $f^{-s}(q)$ for $\frac{1}{f^s(q)}$. 
Also, from now on we will refer to $f^{-*}$ just as the reciprocal of $f$, 
omitting the adjective ``regular'' for the sake of simplicity. 
The reciprocal $f^{-*}$ of a regular function is regular due to
the following result, which can be proven by direct computation.

\begin{lem}
Let $f,g: B = B(0,R) \to \hh$ be regular functions and suppose the
power series expansion of $g$ at $0$, $g(q) = \sum_{n \in \nn} q^n
r_n$, has real coefficients $r_n \in \rr$. Then the function $h : B
\setminus \z_g \to \hh$ defined by $h(q)=\frac{1}{g(q)} f(q)$ is regular.
\end{lem}

If $f$ does not have any zeros then $f^{-*}$ is the inverse element of the
function $f$ with respect to regular multiplication:

\begin{prop}\label{inverse}
Let $f$ be a regular function on $B=B(0,R)$ which does not have
zeros. Then $f^{-*}$ is a regular function on $B$ and
\begin{equation}
f * f^{-*} = f^{-*} * f = 1.
\end{equation}
\end{prop}

\begin{proof} 
By theorem \ref{conjugatezeros}, if $f$ does not have zeros then the
zero set $\z_{f^s}$ of its symmetrization is empty. In this case
$f^{-s}$ and $f^{-*}$ are well defined and regular on the whole ball
$B$ and we may consider their regular product with other regular
functions $g: B \to \hh$. We first note that for all such $g$ we have
$f^{-s}(q) g(q) =  f^{-s}(q)*g(q) = g(q) *  f^{-s}(q)$. We may then
compute: $$f^{-*}*f = f^{-s}* f^c * f =f^{-s} f^s =1,$$  $$f*f^{-*} =
f*f^{-s}*f^c =f^{-s}* f*f^c = f^{-s} f^s =1.$$
\end{proof}

We now give an alternative expression of the regular product and the
reciprocal.

\begin{prop}
Let $f,g$ be regular functions on $B=B(0,R)$. Then
\begin{equation}
f*g(q) = f(q) g\left(f(q)^{-1} q f(q)\right)
\end{equation}
for all $q \in B \setminus \z_{f}$.
\end{prop}

\begin{proof}
Let $f(q) = \sum_{n \in \nn} q^n a_n, g(q) = \sum_{n \in \nn} q^n
b_n$ be the power series expansions of $f$ and $g$. By definition
\ref{multiplication}, $$f*g(q) = \sum_{n \in \nn} q^n \sum_{k=0}^n
a_k b_{n-k} =  \sum_{k,l \in \nn} q^{l+k} a_k b_l =  \sum_{l \in \nn}
q^l f(q) b_l.$$ If $f(q) \neq 0$ then $$\sum_{l \in \nn} q^l f(q) b_l
= \sum_{l \in \nn} f(q) f(q)^{-1} q^l f(q) b_l = f(q) \sum_{l \in
\nn}  \left[f(q)^{-1} q f(q)\right]^l b_l = f(q) g(f(q)^{-1} q
f(q)).$$
\end{proof}

\begin{prop}
Let $f$ be a regular function on $B=B(0,R)$. If we set $T_f(q) =
f^c(q)^{-1} q f^c(q)$, then
\begin{equation}
f^{-*}(q) = f(T_f(q))^{-1}
\end{equation}
for all $q \in B \setminus \z_{f^s}$.
\end{prop}

\begin{proof}
If $f^s(q) \neq 0$ then $f^c(q)\neq 0$, hence $T_f$ is well defined
on $B \setminus \z_{f^s}$. We also note that, by the previous result,
$f^c(q) * g(q) = f^c(q) g(T_f(q))$. We conclude by computation:
$$f^{-*}(q) = f^s(q)^{-1} f^c(q) = [f^c*f(q)]^{-1} f^c(q) =$$ $$=
[f^c(q) f(T_f(q))]^{-1} f^c(q) = f(T_f(q))^{-1} f^c(q)^{-1} f^c(q) =
f(T_f(q))^{-1}.$$
\end{proof}

This result clarifies, together with proposition \ref{inverse}, why 
we call $f^{-*}$ the ``reciprocal'' of $f$.
Remark that $T_f : B \setminus \z_{f^s} \to \hh$ maps any 2-sphere
(or real singleton) $x+y\s$ to itself. In particular $T_f(B \setminus
\z_{f^s}) \subseteq B \setminus \z_{f^s}$. Moreover, we prove the
following.

\begin{prop}
Let $f: B=B(0,R) \to \hh$ be a regular function. $T_f$ and $T_{f^c}$
are mutual inverses. In particular $T_f : B \setminus \z_{f^s} \to B
\setminus \z_{f^s}$ is a diffeomorphism.
\end{prop}

\begin{proof}
For all $q \in B \setminus \z_{f^s}$ we have, setting $p = T_f(q)$,
$$T_{f^c} \circ T_f(q) = T_{f^c}(p) = f(p)^{-1} p f(p) =$$ $$=
f(p)^{-1} \left[f^c(q)^{-1} q f^c(q)\right] f(p) = \left
[f^c(q)f(p)\right]^{-1} q \left [f^c(q)f(p)\right]$$ where
$$f^c(q)f(p) = f^c(q) f(f^c(q)^{-1}q f^c(q)) = f^c* f (q) = f^s
(q).$$ Hence $$T_{f^c} \circ T_f(q) =  f^s(q)^{-1}  q  f^s(q) = q,$$
where the last equality holds because $f^s(q)$ always lies in the
same complex line as $q$ so that they commute.
\end{proof}


\section{The minimum modulus principle}

The reciprocal function will allow us to derive the minimum modulus
principle from the maximum modulus principle proven in
\cite{advances}. Actually, the latter only applies to functions
which are regular on a ball $B(0,R)$. In order to apply it to the
reciprocal $f^{-*} : B \setminus \z_{f^s} \to \hh$ we need
to extend the maximum modulus principle to a slightly larger class of
functions. We first extend the identity principle proven in
\cite{advances}.

\begin{prop}[Identity Principle]\label{identity}
Let $\Omega$ be a domain in $\hh$ intersecting the real axis and
having connected intersection $\Omega_I = \Omega \cap L_I$ with any
complex line $L_I$. If $f,g : \Omega \to \hh$ are regular functions
which coincide on $\Omega \cap \rr$ then they coincide on the whole
domain $\Omega$.
\end{prop}

\begin{proof} Let $h = f-g$ and let us prove $h \equiv 0$ on
$\Omega$. Choose any imaginary unit $I \in \s$ and consider the
restriction $h_I = h_{|_{\Omega_I}}$. Since $h_I : \Omega_I \to \hh$
is holomorphic and it vanishes on the set $\Omega \cap \rr$, which is
not discrete, $h_I$ must vanish identically on $\Omega_I$.
\end{proof}

We are now ready to extend the maximum modulus principle as desired.

\begin{theorem}[Maximum Modulus Principle]\label{maximum}
Let $\Omega$ be the quaternionic ball $B=B(0,R)$ minus a (closed) set
$\mathcal{E}$ consisting of isolated points and isolated 2-spheres $x+y\s$.
Let $f: \Omega \to \hh$ be a regular function. If $|f|$ has a local
maximum point $p \in \Omega$ then $f$ is constant in $\Omega$.
\end{theorem}

\begin{proof} Let $L_I$, with $I \in \s$, be the complex line through
$p$ in $\hh$ and consider the restriction $f_I : \Omega_I \to \hh$.
If $p$ is a maximum point for $|f|$ then it is also a maximum point
for $|f_I|$. By the same argument used in \cite{advances}, we conclude $f_I$ is
constant on $\Omega_I$. By the identity principle \ref{identity}, $f$
is constant on $\Omega$.
\end{proof}

\begin{theorem}[Weak Minimum Modulus Principle]\label{minimum}
Let $f : B=B(0,R) \to \hh$ be a regular function. If $|f|$ has a
local minimum point $p=x+yI \in B$ then either $f$ is constant or
$f$ has a zero in $x+y\s$.
\end{theorem}

\begin{proof}
Suppose $f$ does not have zeroes in $S=x+y\s$ and consider the
reciprocal $f^{-*}$: it is defined on the domain $\Omega = B
\setminus \z_{f^s}$, which includes $S$. Remark that $$|f^{-*}(q)| =
\frac{1}{|f(T_f(q))|}$$ for all $q \in \Omega$. Since $T_f$ is a
diffeomorphism of $B \setminus \z_{f^s}$ onto itself, the fact that
$|f|$ has a minimum at $p = T_f(p')$ implies that $|f \circ T_f|$ has
a minimum at $p'$. As a consequence, $|f^{-*}|$ has a maximum at
$p'$. By the maximum modulus principle \ref{maximum}, $f^{-*}$ is
constant on $\Omega$. This implies that $f$ is constant, too.
\end{proof}

We will soon prove a stronger version of the minimum modulus
principle, which will allow us to prove the open mapping theorem
for regular quaternionic functions. 


\section{Value distribution on the 2-spheres}

We now want to strengthen the minimum modulus principle, proving that
if the modulus of a non-constant regular function $f$ has a minimum at
a non-real point
$p=x+yI$ then $p$ itself must be a zero for $f$. In order to
prove this result, we need to better understand of the distribution
of the values of $f$ on the 2-sphere $x+y\s$ through $p$.

\begin{theorem}\label{values}
Let $f : B=B(0,R) \to \hh$ be a regular function. For all $x,y \in
\rr$ such that $x+y\s \subseteq B$ there exist  $b,c \in
\hh$ such that
\begin{equation}
f(x+yI)=b+I c
\end{equation}
for all $I \in \s$. Remark that $c=0$ if, and only if, $f$ is constant on $x+y\s$.
\end{theorem}

\begin{proof}
Let $f(q)= \sum_{n \in \nn} q^n a_n$ be the power series expansion
of $f$ at $0$ and, for all $n \in \nn$, let $x_n, y_n \in \rr$ be
such that $(x+yI)^n=x_n+y_nI$. Then 
\begin{equation}\label{b e c}
f(x+yI) =  \sum_{n \in \nn}
(x+yI)^n a_n = \sum_{n \in \nn} (x_n + y_n I) a_n =  \sum_{n \in \nn}
x_n a_n + I  \sum_{n \in \nn} y_n a_n =b+Ic
\end{equation}
\end{proof}

We label the 2-spheres characterized by $c=0$ with the following definition.

\begin{definition}\label{degenerate}
Let $f : B=B(0,R) \to \hh$ be a regular function and let $x,y \in
\rr, y>0$ be such that $S=x+y\s \subseteq B$. The 2-sphere $S$ is said to
be \textnormal{degenerate} for $f$ if  the restriction $f_{|_S}$ is
constant. The union $D_{f}$ of all degerate spheres for $f$ is called the
\textnormal{degenerate set} of $f$.
\end{definition}

\begin{ex} \label{esempio1}
    \begin{itemize}
    \item[\textnormal{i)}] \textnormal{Let $P(q)$ be a quaternionic polynomial having real
    coefficients and a non-real root $x+yI$. Then $x+y\s$ is a
    degenerate sphere for $P$. A simple example is given by the
    polynomial $q^{2}+1$, which vanishes on the 2-sphere $\s$. }
    \item[\textnormal{ii)}] \textnormal{The regular function $f: B(0,1) \to \hh$ defined by
    $f(q)=\sum_{n\in \nn}q^{n}= \frac{1}{1-q}$ has no degenerate
    spheres. In fact in this case $f_{I}: B_{I}\to L_{I}$ can be
    written as $f_{I}(x+yI)=\frac{1-x}{(1-x)^{2}+y^{2}}+\frac{y}{(1-x)^{2}+y^{2}}I$
    and $\frac{y}{(1-x)^{2}+y^{2}}=0$ if, and anly if, $y=0$.
    Therefore there are no degenerate spheres for $f$.}
    \end{itemize}
\end{ex}

If $f : B \to \hh$ is regular and $S=x+y\s \subseteq B$ (with $x,y \in
\rr, y>0$) is not degenerate, then the restriction $f_{|_S}$ is an affine
map of $S$ onto a 2-sphere $b+\s c$ with $b,c \in \hh$. As an
immediate consequence, we get the following corollary.

\begin{corollary}
Let $f : B=B(0,R) \to \hh$ be a regular function, let $x,y \in \rr$
be such that $S=x+y\s \subseteq B$ and suppose $S$ not to be
degenerate. Then $|f_{|_S}|$ has one global minimum, one global
maximum and no other extremal point.
\end{corollary}

As a consequence, we get a stronger version of the minimum modulus
principle.
\begin{theorem}[Minimum Modulus Principle]\label{strongminimum}
Let $f : B=B(0,R) \to \hh$ be a regular function. If $|f|$ has a
local minimum point $p\in B$ then either $f(p)=0$ or $f$ is
constant.
\end{theorem}

\begin{proof}
By theorem \ref{minimum} a non-constant $f$ whose modulus has a
minimum at $p=x+yI$ must have a zero at a point $p' \in x+y\s$. In
particular $|f|$ has a minimum at $p'$. By the previous corollary, $|f|$
cannot have two distinct local minimum points on $x+y\s$, unless 
$x+y\s$ is degenerate for $f$. As a consequence, either $f$ is constant on
$x+y\s$ or $p = p'$. In both cases, $f(p) = f(p') = 0$.
\end{proof}

This result will prove very useful in the next section.


\section{The open mapping theorem}\label{opensection}

We are now ready to prove the open mapping theorem for quaternionic
regular functions. Before doing so, we want to
state two properties of the degenerate set of a regular function, defined
in \ref{degenerate}.

\begin{prop}\label{chiusura di Df} If $f : B=B(0, R) \to \hh$ is a
    regular function, then the degenerate set $D_{f}$ of $f$ is a
    closed subset of $B\setminus \rr$. Moreover, if $f$ is not a constant 
    function then the interior of $D_f$ is empty.
\end{prop}
\begin{proof}  By theorem \ref{values}, there exist real analytic
    functions $b=b(x, y), c=c(x,y)$ such that $f(x+yI) = b(x,y)+Ic(x,y)$. 
    Denote by $C_{f}$  the union of all spheres (or real singletons) 
    $x+y\s, x,y \in \rr$ such that $c(x,y)=0$. We have $D_{f} = 
    C_{f}\setminus \rr$ since, for all $x,y \in \rr, y \neq 0$, $f$ is constant on 
    $x+y\s$ if and only if $c(x,y) = 0$. Clearly $C_f$ is closed in $B$, thus 
    $D_{f}=C_{f}\setminus \rr$ is closed in $B\setminus \rr$. Moreover, if 
    $f$ is not constant then the interior of $D_f$ is empty because the interior 
    of $C_f$ is. Assume by contradiction that, for some $I\in \s$, there exist 
    a non-empty open subset $A_I$ of $B_I = B \cap L_I$ such that 
    $A_I \subseteq C_f$, i.e. such that $c(x,y)=  \sum_{n \in \nn} y_n a_n =0$ 
    for all $x+yI\in A_I$. If we set the obvious notation $a_n = 
    a_n^0+a_n^1i+a_n^2j+a_n^3k$ then, for $p=0,1, 2 ,3$, the holomorphic 
    function $f_p: B_I \to L_I$ defined by $$f_p(x+yI) =  
    \sum_{n \in \nn} (x+yI)^n a_n^p = \sum_{n \in \nn} x_n a_n^p + 
    \sum_{n \in \nn} y_n a_n^p I$$ has vanishing imaginary part on the 
    open set $A_I$. Thus $f(A_I) \subseteq \rr$ and, by the open mapping 
    theorem for holomorphic functions of one complex variable, $f_p$ is 
    constant on $B_I$ for  $p=0,1, 2,3$. Hence $f$ is constant on $B_I$ 
    and, by the identity principle \ref{identity}, on $B$.
\end{proof}

\begin{theorem}[Open Mapping Theorem]\label{open}
Let $f : B(0,R) \to \hh$ be a non-constant regular function and let
$D_{f}$ be its degenerate set. Then $f: B(0, R)\setminus \overline{D_{f}}\to
\hh$ is  
open.
\end{theorem}

\begin{proof}
Let $U$ be an open subset of $B(0,R)\setminus \overline{D_{f}} $ and let 
$p_0 \in f(U)$. We will
show that the image $f(U)$ contains a ball $B(p_0,\varepsilon)$ with
$\varepsilon > 0$. Choose $q_0 \in U$ such that $f(q_0) = p_0$. Since
$U\cap D_{f} = \emptyset$, the point $q_0$ must be an isolated zero
of the function $f(q) - p_0$. We may thus choose $r>0$ such that
$\overline{B(q_0,r)} \subseteq U$ and $f(q) - p_0 \neq 0$ for all $q
\in \overline{B(q_0,r)} \setminus \{ q_0 \}$. Let $\varepsilon >0$ be
such that $$|f(q) - p_0| \geq 3 \varepsilon$$ for all $q$ such that
$|q-q_0| = r$. For all $p \in B(p_0, \varepsilon)$ we get $$|f(q) -p|
\geq |f(q) -p_0| - |p-p_0| \geq 3 \varepsilon -\varepsilon = 2
\varepsilon$$ for $|q-q_0| = r$, while $$|f(q_0)-p| = |p_0-p| \leq
\varepsilon.$$ Thus $|f(q_0)-p| < \min_{|q-q_0| = r} |f(q) -p|$ and
$|f(q) -p|$ must have a local minimum point in $B(q_0,r)$. Since
$f(q) -p$ is not constant, it must vanish at the same point by
theorem \ref{strongminimum}. Hence there exists $q \in B(q_0,r)
\subseteq U$ such that $f(q) = p$.
\end{proof}

If $f$ has no degenerate spheres, as in example \ref{esempio1} ii), then 
we get the statement

\begin{corollary} Let $f : B(0,R) \to \hh$ be a regular function. 
If $f$ is not constant on any 2-sphere $x+y\s$ with $y\neq 0$,  then $f$ is  
open.
\end{corollary}

We point out that there exist regular functions $f: B(0, R) \to \hh$ which are
open on $B(0,R)\setminus \overline{D_{f}} $ but not on 
$B(0, R)$. For instance the function
$f(q)=q^{2}+1$ which appears in example \ref{esempio1} i) maps the 
open ball $B=B(I, 1/2)$ centered at $I\in \s$
onto a non-open set: indeed $0\in f(B)$ and $f(B)\cap L_{K}\subseteq \rr$
when $K\in \s$ is orthogonal to $I$. Nevertheless we notice that
$f(B(0,R))$ is open. 
In fact, if we choose an open set $U \subseteq B(0, R)$ which contains
degenerate spheres for $f$, 
we can
still prove that $f(U)$ is open provided $U$  has the following property.

\begin{definition}
Let $U \subseteq \hh$. We say that $U$ is 
\textnormal{circular} if, for all $x+yI \in U$ with $y>0$, the whole
2-sphere $x+y\s$ is contained in $U$.
\end{definition}

For all $x,y \in \rr, r>0$ we may consider the neighborhood of $S =
x+y\s$ defined by $C(S,r) = \{q \in \hh : d(q,S) < r \}$ and remark
that it is circular.

\begin{theorem}\label{circularopen}
Let $f : B(0,R) \to \hh$ be a regular function. If $U$ is a circular open 
subset of $B(0,R)$, then $f(U)$ is open.
\end{theorem}

\begin{proof}
Let $p_0 \in f(U)$. Choose $q_0 = x_0+y_0 I \in U$ such that $f(q_0)
= p_0$. Since $f(q) - p_0$ has a zero on $S=x_0+y_0 \s \subseteq U$,
there exists $r>0$ such that $\overline{C(S,r)} \subseteq U$ and
$f(q) - p_0 \neq 0$ for all $q \in \overline{C(S,r)} \setminus S$.
Let $\varepsilon >0$ be such that $|f(q) - p_0| \geq 3 \varepsilon$
for all $q$ such that $d(q, S) = r$. For all $p \in B(p_0,
\varepsilon)$ we get $$|f(q) -p| \geq |f(q) -p_0| - |p-p_0| \geq 3
\varepsilon -\varepsilon = 2 \varepsilon > \varepsilon \geq |p_0-p| =
|f(q_0)-p| $$ for $d(q, S) = r$. Thus $|f(q) -p|$ must have a local
minimum point in $C(S,r)$. By the minimum modulus principle, there
exists $q \in C(S,r)$ such that $f(q) - p = 0$.
\end{proof}


\end{document}